\documentclass[a4paper,12pt]{amsart}

\usepackage{amsmath}
\usepackage{amssymb, amsbsy, amstext,amsthm}
\usepackage{srcltx}
\usepackage{amscd,amsfonts,color,todonotes,enumerate}
\input amssym.def
\input amssym.tex

\newtheorem{lem}{Lemma}[section]%
\newtheorem{theorem}[lem]{Theorem}%
\newtheorem{cor}[lem]{Corollary}%
\newtheorem{prop}[lem]{Proposition}%

\theoremstyle{definition}
\newtheorem{rem}[lem]{Remark}%
\newtheorem{exam}[lem]{Example}%
\newtheorem{defi}[lem]{Definition}%
\newtheorem{prob}{Problem}%

\def\a{\alpha} \def\b{\beta}  \def\d{\delta} \def\e{\varepsilon}
 \def\s{\sigma} \def\t{\tau}  
 \def\ld{\lambda} 
\def\D{\Delta}  
 \def\O{\Omega} 

\def\oa{\overline A} \def\og{\overline G} \def\oh{\overline H}

\def\on{\overline N}  
  
 \def\ox{\overline X}

   \def\ols{\overline s}
\def\di{\bigm|} \def\lg{\langle} \def\rg{\rangle} \def\lls{\overline \s} \def\olx{\overline x}  \def\oly{\overline y}

\def\PSL{ {\rm PSL}}  
\def\Aut{ {\rm Aut\,}} \def\Inn{ {\rm Inn}}

\def\Mult{ {\rm Mult}}  \def\PGL{ {\rm PGL}}
 \def\AGL{ {\rm AGL}} \def\GL{ {\rm GL}}  \def\PSU{\rm PSU}
 \def\PGammaL{ {\rm P\Gamma L}}  \def\P\GL{ {\rm P\GL}}
 \def\SL{ {\rm SL}} \def\FF{{\mathbb{F}}}
\def\cal{\mathcal} \def\char{ {\rm char}} \def\calm{\mathcal{M}}

\def\ZZ{\mathbb{Z}} 

\def\nd{\mathrel{\bigm|\kern-.7em/}} 
 \def\f{\noindent}
\def\qed{\hfill $\Box$} \def\demo{\f {\bf Proof}\hskip10pt}
\def\gp#1{\langle #1\rangle}

\newenvironment{pf}{ \demo}{\qed \vskip 3mm}

\begin{document}

\title{Skew-morphisms  of  nonabelian characteristically simple groups}

\author{Jiyong Chen}
\address{Southern  University of Science and Technology,  Shenzhen 518055,  P.R.China}
\email{chenjy@sustech.edu.cn}

\author{Shaofei Du}
\address{Capital Normal University, Beijing 100048, P.R.China}
\email{dushf@mail.cnu.edu.cn}

\author{Cai Heng Li}
\address{Southern  University of Science and Technology,  Shenzhen 518055,  P.R.China}
\email{lich@sustech.edu.cn}

\renewcommand{\thefootnote}{\empty}
 \footnotetext{This research was supported by the National Natural Science Foundation of China  under Grant No. 11671276, 11771200, 11931005 and 11971248.}

\begin{abstract}
A skew-morphism  of a finite group $G$ is a permutation $\s$ on $G$ fixing the
identity element, and for which there exists an integer function $\pi$ on $G$ such
that $\s(xy)=\s(x)\s^{\pi(x)}(y)$ for all $x,y\in G$. It has been known that  given a skew-morphism $\s $ of $G$,
 the product of   $\lg \s \rg$  with the left regular
representation of $G$  forms a permutation group on $G$,  called the skew-product group of $\s$. 
The skew-morphism  was  introduced as an algebraic tool to investigate
regular  Cayley maps.   In this paper, the skew-product groups   are characterized, for all  skew-morphisms of   finite nonabelian characteristically simple groups
(see Theorem 1.1) and correspondingly the  Cayley maps on these groups are   characterized (see Theorem 1.5).
 \end{abstract}
\maketitle

\section{Introduction}
Throughout the paper all groups considered are finite, except otherwise stated.

A \textit{skew-morphism} of a group $G$ is a permutation $\s$ of $G$ fixing the
identity element, and for which there exists an integer function $\pi$ on $G$ such that
$\s(gh)=\s(g)\s^{\pi(g)}(h)$ for all $g,h\in G$.

In the special case where
$\pi(g)=1$ for all $g\in G$, the skew-morphism $\s$ is an automorphism of $G$. Thus
skew-morphisms is a generalization of group automorphisms.
Moreover, the investigation of skew-morphisms is at least  related to the following two topics.

\vskip 3mm
(1) {\it  Group factorizations}:    Denote by
$L_{G}:=\{L_{g}\mid g\in G\}$ the left regular representation of $G$.  Then both $\s$
and $L_g$ are permutations on $G$. For any $g,h\in G$,  we have
\begin{equation*}
(\s L_{g})(h)=\s(gh)=\s(g)\s^{\pi(g)}(h)=(L_{\s(g)}\s^{\pi(g)})(h),
\end{equation*}
and so $\s L_{g}=L_{\s(g)}\s^{\pi(g)}$.
Therefore,  $\langle \s\rangle L_{G}\subseteq L_{G}\langle \s\rangle $. Since $G$ is
finite and $L_G\cap \langle\s\rangle =1$,  we have
$\langle \s\rangle L_{G}=L_{G}\langle \s\rangle $, which implies  that $X:=L_{G}\langle\s \rangle$
is a permutation group on $G$, called the \textit{skew-product group} of $\s$, see \cite{CJT2016, ZD2016}.
Since $\langle \s \rangle  $  is a point stabilizer of the transitive
permutation group $X$, it is core-free in $X$, meaning that  any normal subgroup of $X$ contained in $\lg \s\rg $ is trivial.

 Conversely, suppose that  $X$ is any group admitting a complementary factorization $X=GY$
 with $Y=\langle y\rangle$ being cyclic and core-free in $X$. Then for any $g\in G$, there exist unique $g'\in G$
 and $y^{i}\in \langle Y\rangle$ such that $yg=g'y^{i}$.   Define  $\s\colon G\to G$ by
 $\s(g)=g'$, and $\pi\colon G\to \mathbb{Z}$ by $\pi (g)=i$. Then  one may check that $\s$
 is a skew-morphism of $G$ with power function $\pi$.

\vskip 3mm
(2) {\it Cayley maps:}
The concept of skew-morphism was first introduced as an algebraic tool to investigate
regular {\it Cayley maps}~\cite{JS2002}. A Cayley map $\calm =\mathrm{Cay}(G,S,P)$ is an embedding
of a (simple, undirected and connected) Cayley graph of a group $G$ with generating set
$S$ into an orientable closed surface such that, at each vertex $g$ of $\calm $, the local
orientation $R_g$ of the darts $(g,gx)$ incident with $g$ agrees with a prescribed
cyclic permutation $P$ of the generating set $S$, that is, $R_g(g,gx)=(g,gP(x))$ for all
$g\in G$ and $x\in S$. The automorphism group $\Aut(\calm)$ of a Cayley map $\calm$ contains a
vertex-transitive subgroup induced by left multiplication of the elements of $G$. If the
cyclic stabilizer of a vertex is transitive on its adjacent vertices, then the automorphism
group $\Aut(\calm)$ is regular on the darts of $\calm$. In this case the map $\calm$ is called a
\textit{regular Cayley map}. It was shown by Jajcay and \v{S}ir\'{a}\v{n} that a Cayley map $M$ is
regular if and only if $P$ extends to a skew-morphism of $G$. Thus the problem
of determining all regular Cayley maps of a group $G$ is equivalent to the problem of determining
all skew-morphisms of $G$ containing a generating orbit which is closed under taking inverses.

\vskip 3mm
Now we are ready to recall the studying history of skew-morphisms of groups.  An interesting and important problem in this field is a determination of  skew-morphisms
of a given family of finite groups. The problem seems challenging because
even skew-morphisms of the cyclic groups have not yet
been completely determined.  As for the cyclic groups,  skew-morphisms associated
with regular Cayley maps of the cyclic groups can be extracted from the results obtained
in~\cite{CT2014,Kwon2013}. The coset-preserving skew-morphisms of cyclic groups have been completely determined by
Bachrat\'y and Jajcay in~\cite{BJ2014, BJ2016}.
In~\cite{KN2011},
 Kov\'acs and Nedela proved that under certain numerical conditions a skew-morphism of a cyclic
 group may be decomposed into a direct product of skew-morphisms
of cyclic groups of prime power order, then in \cite{KN2017} they determined
the skew-morphisms of cyclic $p$-groups where $p$ is an odd prime, while  that of the case $p=2$ has been recently determined by Hu and the second author \cite{DH}.
 For partial results about skew-morphisms of the dihedral groups,
the reader is referred to~\cite{CJT2016,KK2016,KMM2013,Zhang2015,Zhang20152,ZD2016}. It was told that the corresponding regular Cayley maps have been recently determined by Kan and Kwon.

The skew-morphisms of  finite nonabelian simple groups were determined by  Bachratsy, Conder and  Verret, see \cite{BCV}.
In this paper, we shall extend the results in \cite{BCV} and
 determine    the skew-morphisms of  nonabelian characteristically simple groups.

\vskip 3mm

The main result   of this paper is  the following theorem.

\begin{theorem}\label{main1}
Let $G=T^\ell$ be  a nonabelian  characteristically simple group, where $T$ is simple and $\ell$ is a positive integer. Let
  $X$ be  a skew-product group of a shew morphism $\s$  of $G$. Then  one of following holds:
\begin{itemize}
\item[(1)]  $X=G{:}\lg \s \rg$ and $\s $ is an automorphism of  $G$, or
\item[(2)]  $(X,G)=(\PSL(2,11),A_5)$, $(M_{23},M_{22})$, or $(A_{m+1},A_{m})$ with $m\geqslant 6$ even;  or
\item[(3)]  $X=(T^{\ell-1}{:}\langle \sigma^i\rangle)\times T\lg \tau\rg$ and $\lg\sigma\rg=\lg \sigma^i\rg\times \lg g\tau\rg$, where $(T\lg \tau\rg, T)$ is one of the three pairs in part {\rm (2)}, $i$ is a proper divisor of $|\s|$,  and $g\in T^{\ell-1}$ such that $g^{\s^i}=g$ and $|g|$ is a divisor of $|\tau|$.  Concretely, $X=(A_5^{\ell-1}{:}\lg\s^i\rg)\times \PSL(2,11)$, $(M_{22}^{\ell-1}{:}\lg \s^i\rg)\times M_{23}$, or $(A_m^{\ell-1}{:}\lg\s^i\rg)\times A_{m+1}$.
\end{itemize}
\end{theorem}

\begin{cor}
 Let $G$ be a nonabelian simple group and $X$  a skew-product group of a skew-morphism of $G$. Then either $G$ is normal in $X$ or $(X,G)$ is one of three pairs in (2) of Theorem~\ref{main1}. This coincides with the result in \cite{BCV}.
\end{cor}
 A skew-product $X$ of a group $G$ is called {\it balanced} or {\it simple} if $G\lhd X$ or $X$ is simple, respectively. And $X$ is called {\it mixed} if $X=X_1\times X_2$, where $X_1$ is a balanced skew-product of a skew-morphism of the group $G_1$, $X_2$ is a simple skew-product of a skew-morphism of a group $G_2$ and $G=G_1\times G_2$. Recall that the automorphism group $\Aut(\calm)$ of a regular Cayley map $\calm$ on a group $G$ is a skew-product of a skew-morphism of $G$.
 Correspondingly, a regular Cayley map $\calm$ on a group $G$ is called {\it balanced, simple} or {\it mixed} if the automorphism group $\Aut(\calm)$ is a skew-product of some skew-morphism of the group $G$ which is balanced, simple or mixed, respectively.

  An orientably-regular map $\mathcal M$ can be viewed as a triple $(X,\sigma,\iota)$, where $X$ is the automorphism group of $\mathcal M$ and $(\sigma,\iota)$ is a generating pair for $X$, such that $\lg \sigma\rg$ is a vertex stabilizer of $\mathcal M$ and $\iota$ is an involution interchange the two directions of an edge. See \cite{JS1978} for details.
This map $\mathcal M$ is also denoted by $\mathcal M(X,\sigma,\iota)$.
 Now suppose that $\mathcal M(X,\sigma,\iota)$ is a regular Cayley map on a group $G$. It follows that $X=G\gp{\sigma}$ is a skew-product of the skew-morphism $\sigma$ of $G$. Conversely, if there is an involution $\iota$ in the skew-product $X=G\gp{\sigma}$ such that $\gp{\sigma,\iota}=X$, then there is a regular Cayley map $\calm(X,\sigma,\iota)$ on $G$.

  Now, we apply Theorem~\ref{main1} to characterize regular Cayley maps of characteristically simple group.
  We first introduce some examples which are simple.

 Suppose $\calm$ is a simple regular Cayley map on the group $A_5$. Then, by Theorem~\ref{main1}, $\Aut(\calm)\cong \PSL(2,11)$. Furthermore, there are exactly $5$ non-isomorphic orientably-regular maps, by using the Magma\cite{magma}.

 Suppose $\calm$ is a simple regular Cayley map on the group $M_{22}$. Then, by Theorem~\ref{main1}, $\Aut(\calm)\cong M_{23}$. Furthermore, there are  exactly $330$ non-isomorphic orientably-regular maps, by using the Magma\cite{magma}.

 Suppose $\calm$ is a simple regular Cayley map on the group $A_m$ for an even integer $m\geqslant 6$. Then, by Theorem~\ref{main1}, $\Aut(\calm)\cong A_{m+1}$. Without loss of generality, we may assume that $\sigma=(1,2\dots, m{+}1)$.  There are lots of involutions $\iota\in X$ such that $\lg \sigma, \iota\rg=X$. To enumerate the corresponding maps, a problem naturally arises, which would be independently interesting.

 \begin{prob}
  Characterize and  enumerate involutions $z$ in $A_{m+1}$ such that $\lg (1,2\dots m{+}1), z\rg=A_{m+1}$ for an even integer $m$.
 \end{prob}

 Before showing the whole characterization of regular Cayley maps of characteristically simple groups, we should give following definitions.

\begin{defi}
Let ${\cal M}_1={\cal M}(X_1,\sigma_1,\iota_1)$ and ${\cal M}_2={\cal M}(X_2,\sigma_2,\iota_2)$ be two orientably-regular maps.
The direct product ${\cal M}_1$ and ${\cal M}_2$ is defined to be the map ${\cal M}(X,\sigma,\iota)$ where $\sigma=(\sigma_1,\sigma_2)$, $\iota=(\iota_1,\iota_2)$ are two elements in the direct product group $X_1\times X_2$ and $X=\lg \sigma,\iota\rg\leqslant X_1{\times }X_2$. This map is also denoted by $\calm_1\times \calm_2$. \end{defi}

\begin{rem}
For an orientably-regular map ${\cal M}(X,\sigma,\iota)$, suppose the group $X=X_1\times X_2$ and  let  $\sigma_i$ and $\iota_i$  be the respective projections of $\sigma$ and $\iota$ on subgroup $X_i$ for $i\in \{1,2\}$. Then from $X=\lg \s, \iota\rg\le \lg \sigma_1,\iota_1\rg \times \lg \sigma_2,\iota_2\rg \le X_1\times X_2$  we know that    $X_i=\lg \s_i, \iota _i\rg .$
 Set $\cal M_i=\cal M(X_i, \sigma_i,\iota_i)$. Then $\cal M= \cal M_1\times \cal M_2$. Thus, we say that
that $\cal M$  has a decomposition   $\cal M_1\times \cal M_2$.
\end{rem}

\begin{theorem}\label{thm map}
 Let $G=T^\ell$ where $T$ is a nonabelian simple group and $\ell$ is a positive integer, and let $\cal M$ be a regular Cayley map on $G$. Then $\cal M$ belongs to one of the three cases below.
 \begin{enumerate}
  \item[(1)] $\cal M$ is a balanced Cayley map on $G$;
  \item[(2)] $\cal M$ is a simple Cayley map on $G$, where $G=M_{22}$ or $G=A_m$ where either $m=5$ or $m\geqslant 6$  even;
  \item[(3)] $\cal M$ is a mixed Cayley map on $G$. Furthermore $\calm$ has a decomposition $\cal M=\cal M_1\times \cal M_2$, such that  $\cal M_2$ is one of the maps in case {(2)}, and $\cal M_1$ is an orientably-regular map,  whose automorphism group is a balanced skew-product  but $\cal M_1$   is not necessarily a balanced Cayley map.
 \end{enumerate}
 Moreover, for each case listed above, there are infinitely many regular Cayley maps of this case.
\end{theorem}

 Remind that   an example of a mixed Cayley map will be given  in  Example~\ref{exam-case3}, which cannot be decomposed into a direct product of a balanced Cayley map and a simple Cayley map.

After this introductory section, some preliminary  results  will be given  in Section 2; Theorem~\ref{main1} and Theorem~\ref{thm map} will be  proved in Sections  3 and 4, respectively.

\section{Preliminaries}
In this section we introduce some preliminary results.

 The following  notations will be used in this paper.
The elementary abelian $p$-group of order $p^n$
will be denoted
by $\ZZ_p^n$. Let $q$ be a prime power. Then
the finite field of order $q$ and its corresponding
multiplicative group will be denoted respectively by $\FF_q$ and by $\FF_q^*$.
Let $G$ be a group and $H$ a subgroup of $G.$ Then we use
$G',$  $C_G(H)$ and $N_G(H)$
to denote the derived subgroup of $G$, the centralizer
and the normalizer of $H$ in $G$ respectively.
Let $M$ and $N$ be two groups. Then we use $M\rtimes N$ to denote a
semidirect product of $M$ by $N$, in which $M$ is a normal subgroup.
For a permutation group $G$ on $\O$ and a subset $B$ in $\O$, by $G_B$ and $G_{(B)}$ we denote the setwise stabilizer and pointwise stabilizer, respectively.
\vskip 3mm

A permutation group is called a {\it $c$-group} if it contains  a transitive cyclic regular subgroup.
This is one kind of important  groups, which is related not only to  group theory itself but also to some combinatorial structures.
In \cite{Jon} Jones determined all the primitive $c$-groups and in \cite{LP} Li and Praeger extended it to quasi-primitive groups, almost simple groups and innately simple groups.

\begin{prop}  \label{c-group}  (\cite{Jon}) Every primitive $c$-group is isomorphic to one of the following groups:
$$\ZZ_p\leqslant G\leqslant \AGL(1,p), A_n (n \, {\rm odd}\, ),  S_n, \PGL(d,q)\leqslant  G\leqslant \PGammaL(d,q), M_{11}, \PSL(2,11), M_{23}.$$

\end{prop}

Recall that a  group $G$ is an {\em extension} of $N$ by $H$
if $G$ has a normal subgroup $N$ such that the quotient group
$G/N$ is isomorphic to $H$. In particular, $G$ is a {\em proper
central extension} of $N$  by $H$ if $N \leq Z(G)\cap G'$. For a perfect group $G$  (that is $G=G'$),   such central subgroups are all quotients of a
largest group, called the {\em Schur multiplier} $\Mult(H)$ of $H$.
\vskip 3mm

\begin{prop}  \label{multi}  (\cite{CCN,Suz}) $\Mult (\PSL(2,11))=\ZZ_2$; $\Mult (M_{23})=1$; $\Mult(A_n)=\ZZ_2$ for $5\leqslant n\not\in \{ 6,7\}$, $\ZZ_6$ for $n\in \{6,7\}$.
\end{prop}

\begin{prop}\label{split}
For the simple groups  $\PSL(2,11)$,  $M_{23}$ and $A_{n}$ for $5\leqslant n\neq 6$, we have the following results:
\begin{itemize}
\item[(1)] there exists no nontrivial  proper central extension of $\PSL(2,11)$ containing a subgroup $A_5$;
\item[(2)] there exists no nontrivial  proper central extension of $M_{23}$; and
\item[(3)] there exists no nontrivial  proper  central extension of $A_{n}$ containing a subgroup $A_{n-1}$.
\end{itemize}
\end{prop}
\demo (1) Since $\Mult(\PSL(2,11))=\ZZ_2$ by Proposition~\ref{multi}, there is only one  nontrivial proper central   extension: $\SL(2,11)$. However, there is only one involution in $\SL(2,11)$, so it has no subgroup $A_5$.

\vskip 3mm
(2) Since $\Mult(M_{{23}})=1$, there exists no  nontrivial proper central extension.

\vskip 3mm
(3)   If $n=5$, then $A_5\cong \PSL(2,5)$ and we have a same argument as in (1).


For $n= 7$, the results can be obtained just by Atlas~\cite{CCN}.

Suppose $n\geqslant 8$.
 Let $A$ be a proper central extension of $A_n$. Since   $\Mult (A_n)=\ZZ_2$, we have $A/Z\cong A_n$, where $Z=Z(A)=\lg z\rg\cong\ZZ_2$.
For the contrary,   suppose that $A$ has a subgroup $H\cong A_{n-1}$. Set $\oa=A/Z$ and $\oh=HZ/Z\cong A_{n-1}$.

Since $\oa \cong A_n$ has a faithful right multiplication representation of degree $n$ on $[\oa:\oh]$, we may view that $\oa =A_n$  so that $\oh=A_{n-1}$ as usual.
In fact, in $A_n$, set  $$\overline{s}_1=(1,2,3), \overline{s}_2=(1,2)(3,4), \overline{s}_3=(1,2)(4,5), \cdots, \overline{s}_{n-2}=(1,2)(n{-}1,n).$$
Then $\oa =\lg \overline{s}_1 , \overline{s}_2, \cdots , \overline{s}_{n-2}\rg$ and $\oh =\lg \overline{s}_1 , \overline{s}_2, \cdots , \overline{s}_{n-3}\rg $.
Let $\phi $ be the homomorphism  from $A$ to $A_n$.
Then $\psi=\phi \di_H$ gives an isomorphism from $H$ to $\oh$, while set $s_i=\psi^{-1}(\ols_i)$ for $1\le i\le n-3$.
   Thus $H=\lg s_1,s_2,\dots, s_{n-3}\rg$ with the defining  relations
\[s_1^3=s_i^2=(s_{i-1}s_{i})^3=(s_js_k)^2=1,\]
where $i\geqslant 2, |j-k|\geqslant 2$.

In what follows we shall find an element $s_{n-2}\in A$ such that $\lg H, \s\rg \cong A_n$ so that $A\cong A_n\times \ZZ_2$, achieving a contradiction.
To do this, first, since $(\ols_{n-3}\ols_{n-2})^3=1$ and $z^2=1$, one of two preimages, say $s_{n-2}$  of $\ols_{n-2}$ in $A$  must satisfy   $(s_{n-3}s_{n-2})^3=1$.
Secondly,   for this  $s_{n-2}$ and any $1\leqslant i\leqslant n-4$, take  $x\in A$ such that
\[ (\ols_2\ols_4)^{\olx}=((3,4)(5,6))^{\olx}=(i{+}1,i{+}2)(n{-}1,n)=\overline{s}_i\overline{s}_{n-2},\]
where $\olx=\phi(x)$.  Then from
$(\overline{s}_2\overline{s}_4)^{\overline{x}}=\overline{s}_i\overline{s}_{n-2}$,
 we have   $s_is_{n-2}=(s_2s_4)^xz^{\varepsilon}$, where $\e\in \{0,1\}$  and  then $(s_is_{n-2})^2=((s_2s_4)^2)^x=1$. Finally,
 take  $y\in A$ such that $\oly=\phi(y) =(3,n{-}1)(4,n).$ Then from
   $\overline{s}_2^{\overline{y}}=\ols_{n-2}$ we have  $s_{n-2}=s_2^yz^{\varepsilon'}$ and  then $s_{n-2}^2=(s_2^y)^2=(s_2^2)^y=1$. In summary,
 these elements  $s_i$ ($i=1, 2, \cdots, n-2$)  in $A$ satisfy  exactly  all relations of  $A_n$, that is $\lg s_i\di 1\le i\le n-2\rg \cong A_n$, as desired.
\qed

\section{Proof of Theorem~\ref{main1}}
In this section, we shall prove Theorem~\ref{main1}. To do that, let $G\cong T^\ell$, where $T$ is a nonabelian simple group. Let $X=G\lg \s\rg$ be any skew-product group of $G$.
Let $G_X$ be the core of $G$ in $X$. Then we shall deal with  three cases separately, according to $G_X=G$, $G_X=1$ or $1<G_X<G$ ($\ell\geqslant 2$).
Thus, the proof of Theorem~\ref{main1} consists of the following three lemmas.
\vskip 3mm
\begin{lem} Suppose that $G_X=G.$ Then $G$ satisfies  (1) of Theorem~\ref{main1}.
\end{lem}
 \demo  In this case, $G\lhd X$. Since $G\cong T^\ell$ ($\ell\geqslant 1$), a subgroup $\lg \s^i \rg \lhd X$ if and only if $\s^i\in Z(X)$. Since $\lg \s\rg$ is core-free, we know that the conjugacy action of $\lg \s\rg $ on $G$ is faithful and so $\s$ can be viewed as an automorphism of $G$, as desired.\qed

\begin{lem}\label{2} Suppose that $G_X=1$. Then  $(X,G)$ is given in   (2) of Theorem~\ref{main1}.
\end{lem}
\demo Let $\Omega=[X:G]$ be the set of right cosets of $G$ in $X$. Since $G_X=1$,   $X$ acts faithfully on $\Omega$ by right multiplication.
We shall prove that $X$ must be primitive on $\Omega$. Then $X$ is one   of primitive $c$-groups in   Proposition~\ref{c-group}. Since $X$ is a product of a finite nonabelian characteristically  simple group and a cyclic group,   we know  from Proposition~\ref{c-group} that  $\ell=1$ and  $(X,G)=({\rm PSL}(2,11), A_5)$ or $(M_{23},M_{22})$ or $(A_{m+1},A_{m})$ for some even integer $m\geqslant 6$.

   Without loss of generality, assume that $X$ is the minimal counter example respect to
  $|\O| $,  that is, every  permutation $c$-group of degree less than $|\Omega|$ is primitive, whose   stabilizers are characteristically  simple.
 Let ${\cal B}$ be a non-trivial minimal block system, and let $K$ be the kernel of the action $X$ on $\cal{B}$. Let $B\in {\cal B}$ be  the block containing $w=G$ and set $\lg\sigma\rg _{B}=\lg\sigma^b\rg$, where $b=|{\cal B}|\geqslant 2$. Clearly, $\lg \sigma^b\rg \leqslant K$.  Since  $K_w=X_w\cap K=G\cap K\lhd G$ and $G=T^\ell$ is characteristically  simple, we get $K_w\cong T^{\ell_1}$  where $0\leqslant \ell_1\leqslant \ell$.
  Then $K=K\cap X_{B}=K\cap (G\lg \s^b\rg )=K_w\lg \s^b\rg$. Now  we carry out the proof by the following  five steps:

  \vskip 3mm (i)  Show $K\in \{ \PSL(2,11), M_{23}, A_{m+1}\}$, where $m\geqslant 6$ is  even.
  \vskip 3mm

    Suppose $K_w=1$. Then   $K=\lg \sigma^b\rg \lhd X$, contradicting to $\lg\sigma\rg_X=1$.
So $K_w\ne 1$, that is $\ell_1\geqslant 1$.

   Note that $K=K_w\lg \s^b\rg $ where $K_w\cong T^{\ell_1}$.  If $K_w\lhd K$, then   $K_w \ \char\  K\lhd X$, contradicting with $G_X=1$. Therefore, $K_w\ntriangleleft K$ and so $K^B$   is insolvable, which implies $X_B^B$ are insolvable.  Since $B$ is  a minimal block,  $X_B^B$ is a primitive $c$-group and so $X_B^B\in \{ \PSL(2,11), M_{23}, A_{m+1}\}$, where $m\geqslant 6$ is  even.  Since  $1\ne K^B\lhd X_B^B$, we get  $K^B=X_B^B$.

   Suppose that  $K_{(B)}\ne 1$.
As $K_{(B)}\lhd K_w$, we have $K_{(B)}\cong T^j$ with $1\leqslant j\leqslant \ell_1-1$.
There exists a block $B'\in\cal B$ such that $K_{(B)}^{B'}\ne 1$.
Since $K^{B'}\cong K^B$ is simple and $1\ne K_{(B)}^{B'}\lhd K^{B'}$, we get  $K^{B'}=K_{(B)}^{B'}\cong T^{j'}$ for some $j'\leqslant j$, and it follows from the simplicity of $K^{B'}$  that   $K^{B'}\cong T$ and then $K^B\cong T$.
Since $K_w\cong T^{\ell_1}$, we have $K_w^B=1$  and then $K^B$ is regular on $B$, contradicting to the  insolvability of $K^B$.
Therefore, $K_{(B)}=1$, $\ell_1=1$ and $K\cong K^B\in \{ \PSL(2,11), M_{23}, A_{m+1}\}$, where $m\geqslant 6$ is  even.
 \vskip 3mm
     (ii) Show $X_B=G\lg \s^b\rg =K\times N$, where $N=X_{(B)}\cong T^{\ell-1}$.

    \vskip 3mm
    Set $N=X_{(B)}$. Then  $K\cap N=K_{(B)}=1$. So $KN=K\times N$. Moreover, since $K\cong K^B=X_B^B=X_B/X_{(B)}=X_B/N$, we get $|X_B|=|K||N|$, Therefore,   $X_B=K\times N$ and then $N\cong T^{\ell-1}$, as  $N\lhd G = K_w \times N$.

     \vskip 3mm
     (iii) Show that $\og $  and $\lg \overline \s\rg $ are  core-free in  $\ox=X/K$, in particular, $\ell\geqslant 2$.

     \vskip 3mm
    Since  $\ox$ is a permutation group on ${\cal B}$, with a stabilizer $\ox_B=X_B/K=GK/K=\og$, we know that  $\og$ is core-free in $\ox$.

     To show $\lg \overline \s\rg $ is core-free,  assume the contrary,   $M/K=\lg \overline{\s}\rg _{\ox}\ne \overline 1$.  Then $M\lhd X$ and $M=M\cap (K\lg \s\rg )=K(M\cap \lg \s\rg )$. Set $M\cap \lg \s\rg =\lg \s^c\rg $, where $c$ is a proper factor of $b$, as $\s^b\in K\leqslant M$. Then $M=K\lg \s^c\rg$.

     Now $$\begin{array} {lll}\s^c\in M&\Leftrightarrow &\lg \overline{\s }^c\rg  ^{\og}=\lg \overline{\s }^c\rg \\
     &\Leftrightarrow&  [\overline {\s}^c, \og]=1,  \, \, ({\rm as}\, \, \og\cong T^{\ell-1}) \\
      &\Leftrightarrow&[\lls^c , \on]=1\, \, ({\rm as}\, \, \og=\on) \\
       &\Leftrightarrow& [\s^c,N]\leqslant N\cap K=1\, ({\rm as\,}\, K, N\, \char\, GK).\\
       \end{array}$$  Therefore, since $\lg \s\rg _X=1$ and $[\s^c, N]=1$, any element in $\lg \s^c\rg $  cannot commute with $K$. In other words,   $\lg \s^c \rg $ acts faithfully on $K$ by conjugacy.

   Let $\tau$ be the automorphism of $K$ induced by $\sigma^c$ by conjugation, that is, for any $k\in K$, $k^\tau =k^{\sigma^c}$. Since $\s^c$  is faithful on $K$, we get $|\tau|=|\sigma^c|$. Set $\sigma^b =k_0$, which is  a cyclic regular  subgroup  on $B$, of $K$. Then we get $k^{\tau ^{\frac bc }}=k^{k_0}$ for any $k\in K$. So, $\tau^{\frac bc}=\Inn(k_0)$ and
   $|\tau|=\frac bc |k_0|$ for $\frac bc\geqslant 2$. In what follows, we shall get a contradiction for all three cases of $K$.

   $K=A_{m+1}$ where $m\geqslant 6$ is even:
    since $\Aut(A_{m+1})\cong S_{m+1}$,
     we need to find an element $\tau\in S_{m+1}$ such that $\tau ^{\frac bc}$ is $(m+1)$-cycle $\Inn(k_0)$. Then $\tau $ has to be a $(m+1)$-cycle, and we have  $|\tau|=|\tau ^{\frac bc}|=|k_0|\ne \frac bc|k_0|$, a contradiction.

    $K=\PSL(2,11)$ (resp. $M_{23}$): since $\Aut(\PSL(2, 11))\cong \PGL(2, 11)$ (resp. $M_{23}$), we need to find an element $\tau $ in $\PGL(2,11)$ (resp. $M_{23}$)  such that $|\tau|=11\frac bc$ (resp. $23\frac bc$) for $\frac {b}{c} \geqslant 2$. However, there exists no such element for both groups.

     \vskip 3mm
     (iv) Show the nonexistence of $X$.
   \vskip 3mm
     By (ii) and (iii), we know that $\ell\geqslant 2$, and both  $\og$ and $\lg \overline \s\rg $  are core-free in $\ox$.
     Since $X$ is a minimal counter  example and $|{\cal B}|<|\O|$, we conclude that $X/K=\PSL(2,11)$, $M_{23}$, or $A_{m+1}$ with an even integer $m\geqslant 6$.
Let $L=C_X(K)$. Then $N\leqslant L\lhd X$. Since  $1\ne (K\times L)/K\lhd X/K$, which is simple,  we know that
  $X=K\times L$, where $K,L\in \{PSL(2,11), M_{23}, A_{m+1}\}$, and $\ell=2$. Since $G\cong T^2$, we know that $K\cong L$. However, one cannot find a cyclic regular subgroup of $X$ under the situation.\qed

\begin{lem} Suppose that $1<G_X<G$. Then  $(X,G)$ is given in   (3) of Theorem~\ref{main1}.
\end{lem}
\demo  Set $G_1=G_X\not\in \{ 1, G\}$ and take $G_2$ such that $G=G_1\times G_2$. Then $G_1\cong T^k$, where $1<k<\ell$. Let $\ox=X/G_1$. Suppose that $\ox_1=\lg \overline{\s}\rg _{\ox}=\lg \overline \s^i\rg $ where $i$ is a divisor of $|\s|$. Then $X_1=G_1:\lg \s^i\rg \lhd X$, $\ox=\og_2\lg \overline {\s} \rg$ and $\ox_1\leqslant Z(\ox)$ and  $\ox/\ox_1 =(\og _2 \ox_1 /\ox_1) (\lg\overline \s \rg/\ox_1)$, where both subgroups of $\ox/\ox_1$  are core-free by the definition of $G_1$ and $X_1$.

Using Lemma~\ref{2} to the $c$-group $\ox/\ox_1$, we get that $\ox/\ox_1$ is isomorphic to a simple group  in $\{\PSL(2,11)$, $M_{23}$,  $A_{m+1}$\} with $m\geqslant 6$ even.
  It follows  from $\ox/\ox_1=(\ox/\ox_1)'=\ox' \ox_1/\ox_1$ that $\ox=\ox' \ox_1$, where $\ox'\cap \ox_1$ is a quotient of Schur  multiplier of the corresponding simpe group and $\ox'$ is a central extension of one of three simple groups with the corresponding nonabelian simple subgroup $\og'_2$.
   By Lemma~\ref{split}, the only one possibility is  $\ox'\cap \ox_1=1,$
  that is
 $\ox =\ox' \times \ox_1$, where $\ox'\cong \ox/\ox_1$  must be  isomorphic to  one of $\{ \PSL(2, 11), M_{23}, A_{m+1}\}.$

Since $G_1\leqslant G=G'\leqslant X'$, we have $\ox '=X'/G_1$. Set $X_2=C_{X'}(G_1)$. It follows from $1<(G_2\times G_1)/G_1\leqslant X_2\times G_1/G_1\lhd \ox'$, a simple group,  that $X'=G_1\times X_2$, where $X_2\cong \ox'$, which  is one of $\{ \PSL(2, 11), M_{23}, A_{m+1}\}.$
Then $X=X'X_1=(G_1\times X_2)X_1=X_2X_1$,

 Since $X_2$ is simple not contained in $X_1$, $X_2\cap X_1=1.$ To show $[X_1,X_2]=1$, it  suffices to show $[\s^i,G_2]=1$.
 Recall that $X_1=G_1\lg \s^i\rg$, $X'=G_1\times X_2$ and $G_2\leqslant X_2$. Since $[\ox_1, \ox']=1$,  for any $g_2\in G_2$ we have
   $[\s^i, g_2]\in G_1$.  Since   both $G_1$ and $X_2$ are characteristic subgroups in  $X'$,  $[\s^i, g_2]\leqslant X_2$. This forces
   $[\s^i, g_2]\in G_1\cap X_2=1$. Therefore $[\s^i, G_2]=1$, as desired. Furthermore,  $X=X_1\times X_2$, where $X_1=G_1{:}\lg \s^i\rg$ and
$X_2=X_2\cap ((G_1\times G_2)\lg \s\rg)=G_2(X_2\cap G_1\lg \s\rg )$, which is a simple group in $\{ \PSL(2, 11), M_{23}, A_{m+1}\}$  for $m\geqslant 6$ even, and $G_1\cong G_2^{\ell-1}$.
Set $j=\frac ni$. Then $\lg \s^j\rg\leqslant X'=G_1\times X_2=G\lg \s^j\rg$. Let $\s^j=g\tau$ where $g\in G_1$ and $\tau\in X_2$. It follows that $X_2=X'/G_1=G\lg \s^j\rg/G_1\cong G_2\lg\tau\rg$ and
\[\frac{|X_2|}{|G_2|}\leqslant |\tau|\leqslant |\s^j|=i=\frac{|X_2|}{|G_2|}.\]
This means $|\tau|=i$, $G_2\cap \lg \tau\rg=1$ and $(X_2,G_2)=(T\lg \tau\rg,T)$ is one of the three pairs in part {\rm (2)} of Theorem~\ref{main1}. Furthermore $|g|$ is a divisor of  $|\s^j|=|\tau|$ and $g^{\s^i}=(\s^j\tau^{-1})^{\s^i}=\s^j(\tau^{-1})^{\s^i}=\s^j\tau^{-1}=g$.  \qed

\section{Proof of Theorem~\ref{thm map}}
\demo To prove  Theorem~\ref{thm map},
 let  $G$ be a nonabelian characteristically simple group and $\mathcal M=\cal M(X,\sigma,\iota)$ is an regular Cayley map on $G$.
It follows that $X$ contain a regular subgroup $G$ and $X=G\langle \sigma\rangle$ is a skew-product of $G$.
Therefore, $(X,G)$ belong to one of the three cases stated in Theorem \ref{main1}.

If $X$ is either balanced or simple, then either part (1) or part (2) of Theorem~\ref{thm map} holds.

 Suppose that $X$ is mixed.  Then there exists two subgroups $X_1$, $X_2$ of $X$ such that $X_1\times X_2=X$ where $X_1$ and $X_2$ are skew-products which are balanced and simple respectively. For $i\in \{1,2\}$, let $\s_i$ and $\iota_i$  be the projections of $\s,\iota$ on $X_i$ respectively. Obviously, $\lg \s_i,\iota_i\rg=X_i$. Set $\calm_i=\calm(X_i,\s_i,\iota_i)$. Thus $\calm=\calm_1\times \calm_2$. Moreover,  $\Aut(\calm_1)=X_1$ is a balanced skew-product group and
     $\calm_2$ is  a simple regular Cayley map on $G_2$. Therefore, the part (3) of Theorem~\ref{thm map} holds.

The existences of  infinite many Cayley maps for each of  three cases in Theorem \ref{thm map} will be shown in the following three subsections and thus   the proof of the theorem  is finished.\qed

\subsection{Existence of balanced regular Cayley maps}
The following lemma shows that there are infinitely many balanced regular Cayley maps on characteristically simple groups.
\begin{lem}\label{existence of balanced}
 For any nonabelian simple group $T$ and any  positive integer $\ell$, there exists a regular Cayley map on group $G=T^\ell$ which is balanced.
  \end{lem}
\demo  For any nonabelian simple group $T$ and any positive integer $\ell$, set $G=T^\ell$. For  $\ell\geqslant 2$, we write that $G=T_1\times T_2\times \dots \times T_\ell$ where $T\cong T_i$ by an isomorphism $\varphi_i$ for $1\leqslant i \leqslant \ell$. For  ${\sigma_0}=(1,2, \dots ,\ell)\in S_\ell$, we define its   action on the group $G$  by $\varphi_i(s)^{\sigma_0}=\varphi_{i^{ \sigma_0}}(s)$ for all $s\in T$. Moreover,  let  $\rho_i$  be the projection from $G$ to $T_i$.

In what follows,  we shall show that there exists an automorphism $\sigma \in \Aut(G)$ where $G=T^\ell$   and an involution $\iota\in G$, such that $X=G{:}\lg \sigma\rg=\lg \sigma,\iota\rg$. It follow that $\cal M(X,\sigma,\iota)$ is a balanced regular Cayley map on group $G$.

\vskip 3mm (1) $\ell=1$: Then  $G=T$. By \cite{King}, the nonabelian simple group $T$ can be generated by an involution $\iota$ and an element $r$ of prime order. Set $\sigma=\Inn(r)$ be the inner automorphism induced by $r$, and set $X=G{:}\langle \sigma\rangle.$ Let $N=\langle \iota^{\lg \s\rg}\rg $. Noting that $N$ is normalized by both $\iota$ and $r$, we know that  $1\neq N\lhd G$. Since $G$ is simple, $N=G$ and $X=G{:}\langle \sigma\rangle=\langle \iota,\sigma\rangle.$

\vskip 3mm (2) $\ell=2$:  By   \cite{King} again, let $T=\lg t, r\rg $ where $|t|=2$ and $|r|$ is a prime.    Let $\sigma_1=\Inn({\varphi_1}(r))\in \Aut(G)$ and set
$\sigma=\sigma_1\sigma_0$, $\iota={\varphi_1}(t){\varphi_2}(t)$ and $X=G{:}\lg \sigma\rg $. Let $H=\lg \iota^{\lg \s\rg}\rg .$
Noting that $\lg t^{\lg r\rg}\rg=T$ and $\sigma^2=\Inn(\varphi_1(r) \varphi_2(r))$, the projection map $\rho_1$ and $\rho_2$ are surjective on $H$. By \cite[Lemma, p.328]{scott}, $H$ is either $G$ or a full diagonal subgroup of $G$. Since both $\iota$ and $\iota^\sigma$ are contained in  $H$,   $H$ cannot  be diagonal subgroup and thus $H=G$ so that  $X=\lg \iota,\sigma\rg$.

\vskip 3mm (3) $\ell=3$:  Completely similar arguments to the case when $\ell=2$.

\vskip 3mm (4) $\ell\ge 4$ and $T\not\cong \PSU(3,3)$:  By \cite{MSW1994}, $T$ can be generated by three involutions, namely $x_1,x_2,x_3$. Let $x_{4}=x_{5}=\dots=x_\ell=1\in T$. Now, set $\iota={\varphi_1}(x_1){\varphi_2}(x_2){\varphi_3}(x_3)$, $\sigma=\sigma_0$ and $X=G{:}\lg \sigma\rg$.
  Let $t_{ij}=\varphi_j(x_i)$ $(1\leq i,j\leq \ell)$. One could find that $t_{ij}^\sigma=t_{ij^\sigma}$. Let $H=\langle \iota^{\lg \s\rg}\rg.$  It follows that $H\leq G$ and $H^\sigma=H$. Furthermore
	  \begin{align*}
	    \rho_i(H)&=\rho_i(\langle \iota,\iota^\sigma,\dots, \iota^{\sigma^{\ell-1}}\rangle)\\
	    &=\langle \{\rho_i(\iota^{\sigma^j})|0\leq j\leq \ell-1\}\rangle\\
	    &=\langle \{\rho_i(t_{11^{\sigma^j}}t_{22^{\sigma^j}} t_{33^{\sigma^j}})|0\leq j\leq \ell-1\}\rangle\\
	    &=\langle \{t_{i^{\sigma^{-j}}i}|0\leq j\leq \ell-1\}\rangle\\
	    &=\langle \{{\varphi_i}(x_{i^{\sigma^{-j}}})|0\leq j\leq \ell-1\}\rangle\\
	    &={\varphi_i}(\langle \{x_j|0\leq j\leq \ell-1\}\rangle)\\
	    &=T_i.
	   \end{align*}
	   This means the projection $\rho_i$ is surjective. By \cite[Lemma, p.328]{scott}, $H$ is the direct product $\prod_{j=1}^k H_j$ of full diagonal subgroups  $H_j$ in $\prod_{i\in I_j}T_i$, that is $\rho_i$ is bijective on $H_j$ for each $i\in I_j$,  where the $\{I_1,I_2,\dots , I_k\}$ form a partition of $I=\{1,2,\dots, n\}$. Noting that $H^\sigma=H$,  $\{I_1,I_2,\dots , I_k\}$ form a $\langle \sigma\rangle$-invariant partition. Since $\iota\in H=\prod_{j=1}^k H_j$ is also in the subgroup $\prod_{i=1}^3 T_i$, the set $\{1,2,3 \}$ is a union of some $I_j$. This is impossible if $|I_j|\geq 2$. Hence $|I_j|=1$ and $H=\prod_{i=1}^n T_i=G$. Therefore,  $X=G\langle \sigma\rangle=H\langle \sigma\rangle=\langle \iota,\sigma\rangle$.

\vskip 3mm (5) $\ell\ge 4$ and $T\cong \PSU(3,3)$: If  $\ell=4$, by using the Magma\cite{magma}, there exists an involution $\iota\in G$ such that $\lg \iota, \sigma\rg=X$, where $X=G{:}\lg \sigma\rg$ and $\sigma=\sigma_0$.  Suppose $\ell>4$. Then  by \cite{MSW1994} $T$ can be generated by $4$ involutions, namely $x_1,x_2,x_3,x_4$. Let $x_5=\dots=x_\ell=1\in T$. Now, set $\iota={\varphi_1}(x_1){\varphi_2}(x_2){\varphi_3}(x_3){\varphi_4}(x_4)$, $\sigma=\sigma_0$ and $X=G{:}\lg \sigma\rg$. Using the same argument as in (4), we still have $X=\lg \sigma, \iota\rg$.
\qed

\begin{exam}\label{3psl2p}
 Let $G=\PSL(2,p)$, where $p\geqslant 5$ is a prime.
 Let $\FF_p$ be the finite field of order $p$, $\FF_p^*=\FF_p\setminus\{0\}$ and $\FF'=\{1,2,\dots, \frac{p-1}{2}\}$.
 For any $a\in \FF_p^*$, set $\lambda(a)=\{a,-a\}\cap \FF'$. Let $\D (p)$ be the subset of $\FF'$ consisting of $1$,
	$\ld (\sqrt{2})$ if $p\equiv \pm 1\pmod{8}$, $\ld (\sqrt{3})$ if
	$p\equiv \pm 1\pmod{12},$  and $\ld (\frac{-1\pm \sqrt{5}}2)$ if
	$p\equiv \pm 1\pmod{5}.$   For each $\d \in \FF'\setminus \D (p)$, fix a solution $(\a , \b )$
	to the equation  $\a^2+\b^2-\a+\d \b+1=0.$
	Now, set
	\[ r_\d=\overline{
	\begin{pmatrix}
	\a & \b\\
		\b+\d & 1-\a\\\end{pmatrix}}  \qquad {\rm and } \qquad
	s=\overline{
	\begin{pmatrix}
	0 & 1\\
		{-1} & 0\\\end{pmatrix}}.\]
	By \cite{DW}, $\gp{r_\d, s}=G$ and we have a regular Cayley map $\calm(\PSL(2,p){:}\gp{\Inn(r_\d)},\Inn(r_\d), s)$. Denote this map by $\mathcal{P}(p,\d)$.
\end{exam}

\begin{exam}\label{ppsl2p}
  Let $G=\PSL(2,p)$, where $p\geqslant 5$ is a prime. Set
  \[ r=\overline{\begin{pmatrix}
   1&1\\0&1
  \end{pmatrix}} \qquad \mbox{and}\qquad s_c=\overline{\begin{pmatrix}
  0&-\frac{1}{c} \\c&0\end{pmatrix}},\]
  where $1\leqslant c\leqslant \frac{p-1}{2}$. Note that the only maximal subgroup in $G$ which contain $r$ is $N_G(\gp{r})$. Since $s_c\notin N_G(\gp{r})$, $\gp{r,s_c}=G$. Furthermore, we have a regular Cayley map $\calm(\PSL(2,p){:}\Inn(r),\Inn(r), s_c)$. Denote this map by $\mathcal{Q}(p,c)$.
\end{exam}

\begin{lem}
 Let $G=\PSL(2,p)$, where $p\geqslant 5$ is a prime. Suppose that $\calm$ is a balanced regular Cayley map on $G$. If the valancy of $\calm$ is $3$, then $\calm\cong \mathcal{P}(p,\delta)$ for some $\d$ which is defined in Example~\ref{3psl2p}. If the valancy of $\calm$ is $p$, then $\calm\cong \mathcal{Q}(p,c)$ for some $c$ which is defined in Example~\ref{ppsl2p}.
\end{lem}

\begin{pf}
 Suppose that $\calm=\calm(G{:} \gp{\s},\s,\iota)$ and $|\s|=3$ or $p$. Since
   ${\rm Out}(G)\cong \mathbb{Z}_2$ and $|\s|$ is odd,  we may  view  $\s$ as  an element of $\Inn(G)$. Clearly,  $\iota\in G$. Suppose $\s=\Inn(r)$ and $s=\iota$. Then $r,s\in G$ and $\lg r,s\rg=G$. Furthermore $|r|=|\s|$ is the valency of $\calm$.

	\vskip 3mm (1) $|r|=3$:

	\vskip 3mm
	In \cite{DW}, the set of the representatives of the orbits of $(2,3)$-generating pairs of $\PSL(2,p)$ under of the action of $\PGL(2,p)$ is classified, that is $\{(r_\d, s)\di \d \in \FF'\setminus \D (p)\}$ which is defined in Example~\ref{3psl2p}.
	Hence, we choose $\s =\gp{\Inn(r_\d)}$ and $\iota=s$, so that   $\calm\cong \calm(G\rtimes\gp{\Inn(r_\d)}, \Inn(r_\d), s)\cong \mathcal{P}(p,\delta)$.

	\vskip 3mm
	(2) $|r|=p$.
	\vskip 3mm

	Since there is only one conjugacy class of elements of order $p$ in $\PGL(2,p)$, we set $r=\overline{
	\begin{pmatrix}
	 1 & 1\\
	 	0 & 1\\\end{pmatrix} }$ so that $\s=\Inn(r)$.  The $C=C_{\PGL(2,p)}(r)=\lg r\rg $. Since there is only one class of involutions in $G$, we take $e=\overline{
	\begin{pmatrix}
	0& 1\\
		-1 & 0\\\end{pmatrix}} $ so that $s=e^g$ for some $g\in G$.

	Suppose that $p\equiv 1\pmod{4}$. Then $C_G(e)=D_{p-1}$ and so there are $\frac {p(p+1)}2$ involutions but $p$ of them are  contained in a maximal subgroup $N_G(\lg r\rg)$ which is the unique maximal subgroup of $G$ containing $r$, so there are  $\frac {p(p-1)}2$ of them remaining,  which in turn  form $\frac {p-1}2$ $C$-conjugacy orbits
	Thus we have  $\frac {p-1}2$  choices for $s$.

	Suppose that $p\equiv 3\pmod{4}$. Then $C_G(e)=D_{p+1}$ and so there are $\frac {p(p-1)}2$ involutions,  which  forms $\frac {p-1}2$ $C$-conjugacy orbits.
	Thus again we have  $\frac {p-1}2$  choices for $s$.

	Now we are  finding  such involution $s$. Every involution in $G$  is of form $s=\overline{
	\begin{pmatrix}
	a & b\\
		c & -a\\
		\end{pmatrix}} ,$ where $a^2+bc=-1$. If $c=0$, then   $p\equiv 1\pmod{4}$  and   $s\in N_G(\lg r\rg )$. So we set $c\ne 0$ for all $p$. Now for any $r^t\in C$,
	$s^{r^t}=\overline{\begin{pmatrix}
		a-tc & 2at-ct^2-\frac {a^2+1}{c}\\
		c & tc-a\\
		\end{pmatrix}} .$ For each $c\ne 0$, choose $t=\frac ac$. Then $s_c:=s^{r^{\frac ac}}=\overline{\begin{pmatrix}
			0 & -\frac 1c\\
			c & 0\\\end{pmatrix} } ,$  where $c$ and $-c$ determine the same element in $G$. Therefore, under the action of $\PGL(2,p)$, the set of representatives of $(p,2)$-generating pairs is
 $$\{ r, s_c\di  1\leqslant c\leqslant \frac{p-1}2\}.$$
	Therefore, we take $\iota=s_c$ so that  the  regular Cayley maps $\calm\cong\mathcal{Q}(p,c)$ for some $1\leqslant c\leqslant \frac{p-1}2$.
\end{pf}

\subsection{Existence of simple regular Cayley maps}

In this subsection, we consider the regular Cayley maps with the automorphism group $\PSL(2,11)$, $M_{23}$ and $A_{m+1}$.

\begin{lem}
 There exist exactly  $5$ and $330$ non-isomorphic simple  regular Cayley maps on $A_5$ and $M_{22}$, respectively.
 Furthermore, the face valency  of these maps are listed in Table~\ref{table for psl211m23}.
\end{lem}

\begin{table}[ht]
\begin{center}
 \begin{tabular}{|c|c|c|}
 \hline
 & \# non-isomorphic maps& face valency\\\hline\hline
 $A_5$&5&$3,5^{(2)},16,11$\\\hline
 $M_{22}$&330&$4^{(4)}, 5^{(12)}, 6^{(28)}, 7^{(32)}, 8^{(32)}, 11^{(68)}, 14^{(64)}, 15^{(56)}, 23^{(34)}$\\\hline

 \end{tabular}\end{center}
 \caption{Simple regular Cayley maps on $A_5$ and $M_{22}$}
 \label{table for psl211m23}
\end{table}
The notation $a^{(b)}$ in Table \ref{table for psl211m23} means that there are $b$ non-isomorphic maps with face valency $a$.

\begin{pf}
 Let $\calm =\calm(X,\sigma,\iota)$ be a simple regular Cayley map on $G\in \{A_5, M_{22}\}$. Then $X=G\gp{\s}$ is a simple skew-product of $G$ and $\gp{\s,\iota}=X$. By Theorem~\ref{main1}, $X=\PSL(2,11)$ and $M_{23}$ respectively. Set \[\Delta=\{(\s,\iota)|X=G\gp{\s}=\gp{\s,\iota}, \gp{\s}\cap G=1, |\iota|=2\}.\] The number $n$ of non-isomorphic simple regular Cayley maps on $G$ equals to the number of orbits of $\Aut(X)$ acting on $\Delta$. Noting that $\Aut(X)$ acts on the generating pairs semiregularly, $n=\frac{|\Delta|}{|\Aut(X)|}$.

 For $G=A_5$, we have $X=\PSL(2,11)$ and $|\s|=11$. By checking \cite{CCN}, there are two conjugacy classes of elements with order $11$ and one conjugacy class of involutions.
 Since there are no maximal subgroup of order mulitiple of $22$, all $(11,2)$ pairs generate the group $X$.
 Note that the order of the centralizers of an involution and an element of order $11$ is $12$ and $11$ respectively. We have $|\Delta|=\frac 1{12}|X|\cdot \frac{2}{11}|X|=\frac{1}{66}|X|^2$.
 Thus $n=\frac{|\Delta|}{|\Aut(X)|}=5$.

For $G=M_{22}$, we have $X=M_{23}$ and $|\s|=23$. By similar argument as above, we have $|\Delta|=\frac 1{2688}|X|\cdot \frac{2}{23}|X|=\frac{1}{1344\cdot 23}|X|^2$ and $n=\frac{|\Delta|}{|\Aut(X)|}=\frac{|X|}{1344\cdot 23}=330$.

Finally, the face valency of these maps can be computed by  using the Magma \cite{magma}.
\end{pf}

Now, we consider simple regular Cayley maps on $G=A_m$ for $m\geqslant 6$ an even integer. By Theorem~\ref{main1}, $X=A_{m+1}$ and $|\s|=m+1$ is a full cycle in $X$. Without loss of generality, we assume that $\s=(1,2,3,\dots m{+}1)$. Then there exists a simple regular Cayley map on $G$ if and only if there is an involution $\iota$ in $X$ such that $\gp{\s,\iota}=X$.

\begin{lem}\label{am+1}
 Let $X=A_{m+1}$ for an even integer $m\geqslant 6$ with a natural action on set $\Omega=\{1,2,\dots,m+1\}$. Suppose that $\sigma= (1,2,3,\dots m+1)$ and $\iota=(1,2)(3,4)$ then $X=\gp{\s,\iota}$.
\end{lem}

\begin{pf}
 Let $H_i$ be the subgroups of $X$ that fix $i+1,i+2,\dots,m+1$ where $i=4,5,\dots, m$ and set $H_{m+1}=X$. Obviously, $H_i\cong A_i$. Let $a=\iota^{\sigma}\iota^{\sigma^3}\iota^{\sigma^5}\dots \iota^{\sigma^{m-1}}$. By a simple calculation, we have $a=(1,2,3)$, and $\lg a, \iota\rg=H_4\leqslant \lg \sigma,\iota\rg$. Note that $\lg H_{i}, \iota^{\sigma^{i-3}}\rg=H_{i+1}\leqslant \lg \sigma,\iota\rg$. It follows that  $X=\lg \sigma,\iota \rg$.
\end{pf}

By this lemma, for each even integer $m\geqslant 6$, there is a simple regular Cayley map $\calm(A_{m+1}, \sigma,\iota)$ on group $G=A_{m}$. We have the following corollary.
\begin{cor}
 There are infinitely many simple regular Cayley maps.
\end{cor}

\subsection{Existence of mixed regular Cayley maps}

Since a mixed skew-product group is a direct product of a balanced one and a simple one, it is natural to give examples of mixed regular Cayley maps as the direct product maps.

\begin{lem}\label{existence of mixed}
 Let $m\geqslant 6, \ell\geqslant 5$ be two integers such that $m$ is even and $\gcd(\ell-1, m+1)=1$. Then there is a mixed regular Cayley map on the group $A_m^\ell$. It follows that there are infinitely many mixed regular Cayley maps on characteristically simple group.
\end{lem}

\begin{pf} Set $G_1=A_m^{\ell-1}$. Then by the proof of Lemma~\ref{existence of balanced}, there exist $\sigma\in \Aut(G_1)$ and $\iota_1\in G_1$ such that $|\s_1|=\ell-1$, $|\iota_1|=2$ and  $\gp{\s_1,\iota_2}=G_1{:}\gp{\s_1}$. Set $X_1=G_1{:}\s_1$.

 Let $G_2=A_m<X_2=A_{m+1}$ and let $\s_2=(1,2,\dots,m+1 )$, $\iota_2=(1,2)(2,4)$ be two elements in $X_2$. By Lemma~\ref{am+1}, $X_2=\gp{\s_2,\iota_2}=G_2\gp{\s_2}$, where $G_2\cap\gp{\s_2}=1$.

 Now, Set $X=X_1\times X_2$ and $A_m^\ell \cong G=G_1\times G_2<X$. Let $\s=\s_1\s_2$ and $\iota=\iota_1\iota_2$ be two elements of $X$. Since $\gcd(\ell-1, m+1)=1$, $\gp{\s}=\gp{\s_1}\times \gp{\s_2}$. Note that $G\cap \gp{\s}=1$ and \[X=G_1{:}\gp{\s_1}\times G_2\gp{\s_2}=(G_1\times G_2)(\gp{\s_1}\times \gp{\s_2})=G\gp{\s}. \]
 If follows that $X$ is a mixed skew-product group of $G$ with the skew morphism $\s$.

 Recalling that $\gp{\s}=\gp{\s_1}\times \gp{\s_2}$, we have $\s_1,\s_2\in \gp{\s,\iota}$. Then $\gp{\s_2,\s_2^\iota}=X_2=\gp{\s_2,\iota_2}\leqslant \gp{\s,\iota}$. Hence $\iota_2\in \gp{\s,\iota}$ and $\iota_1=\iota\iota_2\in \gp{\s,\iota}$. It follows that \[\gp{\s,\iota}\geqslant \gp{\s_1,\s_2,\iota_1,\iota_2}=X\geqslant \gp{\s,\iota},\]
 This means $X=\gp{\s,\iota}$.
 Therefore, there is a map $\calm=\calm(X,\s,\iota)=\calm(X_1, \s_1,\iota_1)\times \cal(X_2,\s_2,\iota_2)$ is a mixed regular Cayley map on $G=A_m^\ell$.
\end{pf}

The mixed regular Cayley maps constructed in the proof of Lemma~\ref{existence of mixed} are direct products of simple and balanced regular Cayley maps. This is not always true. There are infinitely many mixed  regular Cayley maps which can not be a direct product of simple and balanced regular Cayley maps, as we mentioned in Theorem \ref{thm map}.

\begin{exam}\label{exam-case3}
 Let $n\geqslant 3$ be an odd integer. By the Dirichlet prime number theorem, there are infinitely many prime $p$ such that $n<p$ and $p\equiv -1\pmod{n}$.
 Let $m=n+p$, $\O=\{1,2,\dots,m\}$  and $\overline{\O}=\{\bar{1},\bar{2},\dots,\overline{m+1}\}$. Pick up  the following elements from the alternating group $A_\O$ and $A_{\overline {\O}}$:
 \begin{align*}
 r_1&=(1,2,\dots ,p)\in A_\O;\\
 s_1&=(1,p+1)(2,p+2)\dots (n,p+n)(n+1 ,n+2)\in A_\O;\\
 g_1&=(p+1,p+2,\dots,p+n)\in A_\O;\\
 r_2&=(\bar{1},\bar{2},\dots,\overline{n+p+1})\in A_{\overline {\O}};\\
 s_2&=(\bar{1},\bar{2})(\bar{3},\bar{4})\in A_{\overline {\O}}.
 \end{align*}
 Now set $G=A_\O\times A_{\overline {\O}\setminus\{\bar{1}\}}\cong A_m^2$ and
 \begin{align*}
  X_1& =A_\O{:}\gp{\Inn(r_1)},&
  \s_1& =\Inn(r_1)g_1\in X_1,&
  \iota_1& =s_1\in X_1,\\
  X_2& =A_{\overline {\O}},&
  \s_2& =r_2 \in X_2,&
  \iota_2& =s_2\in X_2,\\
  X& =X_1\times X_2,&
  \s& =\s_1\s_2\in X,&
  \iota& =\iota_1\iota_2\in X.
 \end{align*}
Clearly, $X_1$ is a balanced skew-product group and $X_2$ is a simple one.
We shall show that $X=\gp{\s,\iota}$ so that  the map $\calm=\calm(X, \s, \iota)$ is a mixed regular Cayley map on $G$, and that
   $\calm$ cannot be decomposed into  a direct product  of  a balanced Cayley map and a  simple Cayley map.

We divide the proof of our claim into the following three steps:

\vskip 3mm (1) The group $\gp{r_1,s_1}=A_\O$ and $X_1=\gp{\s_1, \iota_1}$.
\vskip 3mm

 Obviously,  $\gp{r_1,s_1}$ is transitive on $\O$. Following  \cite[Theorem 13.9]{wielandt-book},  to show $\gp{r_1,s_1}=A_\O$, it is sufficed to show that  $\gp{r_1,s_1}$ acts primitively   on $\O$.  In fact, take any block  $\Delta$,  containing 1.  Then we have either  $\Delta^{r_1}=\Delta$ or $\Delta^{r_1}\cap \Delta=\emptyset$, which implies  $\{ 1, 2, \cdots , p\}$ is contained either in $\Delta$ or in $p$ distinct blocks. The first case implies $|\Delta|>\frac 12|\Omega|$ and
so $\D =\O $, and the second case implies   $|\Delta|\leqslant \frac{|\O|}{p}<2$, meaning
  $\Delta=\{1\}$. In both cases, $\D$ is just trivial.

 Since $[\Inn(r_1),g_1]=1$ and $\gcd(|\Inn(r_1)|,|g_1|)=1$, we have  $\gp{\s_1}=\lg \Inn(r_1)\rg \times \lg g_1\rg$. Furthermore,
 \begin{align*}
 X_1&\geqslant  \gp{\s_1,\iota_1}\geqslant \gp{\Inn{(r_1)},s_1}\geqslant \gp{s_1^{\gp{r_1}}}{:}\gp{\Inn{(r_1)}}=A_\O{:}\gp{\Inn(r_1)}=X_1. \end{align*} This means $X_1=\gp{\s_1, \iota_1}$.

 \vskip 3mm (2) The group $X=\gp{\s,\iota}=G\gp{\s}$ is a mixed skew-product of $G$.
 \vskip 3mm

  By Lemma~\ref{am+1}, $\gp{\s_2, \iota_2}=X_2=A_{\overline{\O}\setminus\{\bar{1}\}}\gp{\s_2}$ is a simple skew-product group. Let $\rho_1$ be the projection from $X=X_1\times X_2$ to $X_1$ and let $\rho$ be the restriction of $\rho_1$ on the subgroup $\gp{\s,\iota}$. Then the image ${\rm Im}\rho=\gp{\s_1,\iota_1}=X_1$ and the kernel $\ker \rho\lhd \gp{\s,\iota}$ . Since $\ker\rho\leqslant X_2$ and both $\s_1$ and $\iota_1$ are commute with  $\ker\rho$,   both  $\s\s_1^{-1}$ and $\iota \iota_1$ normalize $\ker\rho$,  which implies $\ker\rho \lhd X_2$. Since $|\s|\ge |\s_1|$, we know that $\ker\rho \ne 1$ but   $\ker\rho=X_2$. Thus, $\gp{\s,\iota}=X$.
  It is easy to verify that $G\cap \gp{\s}=1$. Comparing the order of $X$ and $G\gp{\s}$, we have $X=G\gp{\s}$. Moreover, $X$ is a mixed skew-product of the skew-morphism $\s$.

 \vskip 3mm (3) $\calm$ can not be a direct product  of two regular Cayley maps which are simple and balanced respectively.
\vskip 3mm

  For the contrary, suppose that $\calm(X,\s,\iota)=\calm(X_3,\s_3,\iota_3)\times \calm(X_4,\s_4, \iota_4) $ where $\calm(X_3,\s_3,\iota_3)$ is balanced and $\calm(X_4,\s_4,\iota_4)$ is simple. By a similar argument in (2), $X=X_3\times X_4=\gp{\s_3\s_4, \iota_3\iota_4}$.
  Thus there is a group automorphism $\tau$ of $X$ such that $\s^\tau=\s_3\s_4$.
      Note that  $X$ has only two  nonabelian minimal normal subgroups: $A_{\O}$ and  $A_{\overline{\O}},$ while both are preserved by $\t$.  Since $X'=X_3'\times X_4'=X_3'\times X_4$, we know that
      $X_3'=A_{\O}$ and $X_4=A_{\overline{\O}}.$
       Moreover,  $\frac{|X|}{|X'|}=\frac{|X_1|}{|X_1'|}=|r_1|=p=\frac{|X_3|}{X_3'}=|\s_3|$ and  $\s^p=(\Inn(r_1)g_1\s_2)^p=g_1^p\s_2^p=g_1^pr_2^p$ is contained
         neither in $A_{\O}$ nor in $ A_{\overline{\O}}$. But
   $(\s^p)^\tau=\s_3^p\s_4^p=\s_4^p\in X_4=A_{\overline{\O}}$, that is $\s^p\in A_{\overline{\O}}^{\t^{-1}}=A_{\overline{\O}}$, a contradiction. \qed
\end{exam}


\end{document}